\documentclass[12pt]{amsart}
\begin{document}
\def\R{{\mathbb R}}
\def\Rad{{\mathcal R}}
\def\Z{{\mathbb Z}}
\def\C{{\mathbb C}}
\newcommand{\trace}{\rm trace}
\newcommand{\Ex}{{\mathbb{E}}}
\newcommand{\Prob}{{\mathbb{P}}}
\newcommand{\E}{{\cal E}}
\newcommand{\F}{{\cal F}}
\newtheorem{df}{Definition}
\newtheorem{theorem}{Theorem}
\newtheorem{lemma}{Lemma}        
\newtheorem{pr}{Proposition}
\newtheorem{co}{Corollary}
\newtheorem{pb}{Problem}
\newtheorem{ex}{Example}
\newtheorem{remark}{Remark}
\def\n{\nu}
\def\sign{\mbox{ sign }}
\def\a{\alpha}
\def\N{{\mathbb N}}
\def\A{{\cal A}}
\def\L{{\cal L}}
\def\X{{\cal X}}
\def\F{{\cal F}}
\def\c{\bar{c}}
\def\v{\nu}
\def\d{\delta}
\def\diam{\mbox{\rm dim}}
\def\vol{\mbox{\rm Vol}}
\def\b{\beta}
\def\t{\theta}
\def\l{\lambda}
\def\e{\varepsilon}
\def\colon{{:}\;}
\def\pf{\noindent {\bf Proof :  \  }}
\def\endpf{ \begin{flushright}
$ \Box $ \\
\end{flushright}}
\title[Functions positively associated with integral transforms]
{Functions positively associated with integral transforms}

\author{Alexander Koldobsky}
\address {Department of Mathematics, University of Missouri, Columbia, MO 65211, USA}
\email{koldobskiya@@missouri.edu}

\thanks{The author was supported in part by the U.S. National Science Foundation Grant DMS-2450745.}
\date{}

\begin{abstract}  We introduce the class of functions positively associated with a linear operator. We describe these classes for several integral operators including the $q$-cosine transform and the spherical Radon transform. We show that positively associated functions control the comparison problem for linear operators generalizing the Busemann-Petty problem for convex bodies.
\end{abstract}  
\maketitle
\section{Introduction}
In this note we consider the following problem.
\begin{pb}\label{general} Let $T: C_e^\infty(S^{n-1})\to C_e^\infty(S^{n-1})$ be a linear operator, where $C_e^\infty(S^{n-1})$ is the space of even real valued infinitely differentiable functions on the  sphere $S^{n-1}$ in $\R^n.$ Suppose that $f,g\in C_e^\infty(S^{n-1})$ are positive functions, and $Tf(x)\le Tg(x)$ for every $x\in S^{n-1}.$ Does this condition alone allow to compare the $L_p$-norms  of the functions $f$ and $g$ for a given $p>1$?
\end{pb}
\smallbreak
The motivation comes from a connection with convex geometry. The Busemann-Petty problem posed in \cite{BP} in 1956 asks whether origin-symmetric convex bodies in $\R^n$ with uniformly smaller areas of central hyperplane sections necessarily have smaller volume. More precisely, suppose that $K$ and $L$ are origin-symmetric convex bodies in $\R^n,$ and 
$|K\cap \xi^\bot|\le |L\cap \xi^\bot|$ for every $\xi\in S^{n-1}$, where $\xi^\bot$ is the central hyperplane perpendicular to $\xi.$
Does it follow that $|K|\le |L|$? Here $|K|$ stands for volume in the appropriate dimension. The problem was solved at the end of the 1990's, as the result of work of many mathematicians. The answer is affirmative if the dimension $n\le 4,$ and it is negative if $n\ge 5;$ see \cite{G3,K1} for the solution and its history.

Lutwak \cite{Lu} introduced the class ${\mathcal I}_n$ of intersection bodies that controls the Busemann-Petty problem in the sense that
if $K\in {\mathcal I}_n$ then the answer is affirmative for any body $L$ with greater central hyperplane sections. On the other hand, if $L\notin {\mathcal I}_n$ then one can perturb $L$ to construct a counterexample. Lutwak's connection reduced the problem to the study of intersection bodies, and the final solution follows from the fact that every origin-symmetric convex body in $\R^n,\ n\le 4$ is an intersection body, while in dimensions 5 and higher there
exist origin-symmetric convex bodies that are not intersection bodies; see \cite{Ga2, GKS, Z3}.

The Busemann-Petty problem is a particular case of Problem \ref{general}. To see this, consider the case where $T=R$ is the spherical Radon transform defined by
\begin{equation}\label{spherical-Radon}
Rf(\xi)= \int_{S^{n-1}\cap \xi^\bot} f(x) dx,\quad \forall \xi\in S^{n-1},\quad \forall f\in C(S^{n-1}),
\end{equation}
$p=\frac n{n-1},$ and $f=\|\cdot\|_K^{-n+1},\ g=\|\cdot\|_L^{-n+1}.$ Then the question of Problem 1 turns into the Busemann-Petty problem. Note that in the solution of the Busemann-Petty problem in \cite{GKS} most of the work is done for infinitely smooth bodies, and the general result follows by approximation.
\smallbreak
We will show that there exists a class ${\rm Pos}(T)$ of functions positively associated with the operator $T$ that controls Problem \ref{general} in the same way as intersection bodies control the Busemann-Petty problem. Namely, if $f^{p-1}\in {\rm Pos}(T)$, then the $L_p$-norm of $f$ is smaller than the $L_p$-norm of $g$ provided that $Tf\le Tg$ pointwise. On the other hand, if $g^{p-1}\notin {\rm Pos}(T)$, then one can construct a function $f$ giving together with $g$ a counterexample. We will describe ${\rm Pos}(T)$ for a class of integral operators including $q$-cosine transforms 
\begin{equation}\label{q-cosine}
C_qf(\xi)=\int_{S^{n-1}} |(x,\xi)|^q f(x) dx.
\end{equation}
We get the following answer to Problem \ref{general} for $q$-cosine transforms.
\begin{theorem}\label{comp-q>0}  Let $p>1,\ q>0,\ q\neq 2,4,6,...$ 
\smallbreak
(i) Suppose that $f\in C_e^\infty(S^{n-1})$ is a positive function such that the Fourier transform of $\frac 1{\Gamma(-q/2)} f^{p-1}\cdot r^q$ is a positive distribution on $\R^n\setminus \{0\}$. If $g\in C_e^\infty(S^{n-1})$ is a positive function such that $C_qf\le C_qg$ pointwise on $S^{n-1}$, then $$\|f\|_{L_p(S^{n-1})} \le \|g\|_{L_p(S^{n-1})}.$$

(ii) If $g\in C^\infty_e(S^{n-1})$ is a strictly positive function such that the Fourier transform of $\frac 1{\Gamma(-q/2)}g^{p-1}\cdot r^q$ is not a positive distribution on $\R^n\setminus \{0\}$, then there exists a strictly positive function $f\in C^\infty_e(S^{n-1})$ such that $C_qf\le C_qg$ pointwise on $S^{n-1}$, but
$$\|f\|_{L_p(S^{n-1})} > \|g\|_{L_p(S^{n-1})}.$$
\end{theorem}
The condition of Theorem \ref{comp-q>0} admits a geometric interpretation. For $q>0,\ q\neq 2,4,6...$ and a strictly positive function $f$, the Fourier transform of $\frac 1{\Gamma(-q/2)} f^{p-1}\cdot r^q$ is a positive distribution on $\R^n\setminus \{0\}$ if and only if $f=\|\cdot\|^{\frac q{p-1}}$ where $(\R^n,\|\cdot\|)$ is an $n$-dimensional subspace of $L_q;$ see \cite{K-1992}.

\smallbreak

\section{Positively associated functions}

\begin{df} \label{positive} Let $X,Y$ be ordered topological vector spaces, let $X^{'}$ be the dual space of continuous linear functionals on $X,$ and let $H$ be a subspace of $X^{'}.$ Let $T: X\to Y$ be a linear operator. 
We say that a functional $h\in H$ is {\bf positively associated} with the operator $T$ if for any $\phi \in X$ satisfying  $T\phi \ge 0$ we have $\langle h, \phi \rangle \ge 0.$
\end{df}

We denote the set of such functionals $h$ by ${\rm Pos}(T)={\rm Pos}(T; X,Y,H).$ The latter equality means that the definition includes the choice of $X,Y,H.$

We start with an example. Let $X=Y={\mathcal S}_e={\mathcal S}_e(\R^n)$ be the space of real valued even functions $\phi\in C^\infty(\R^n)$ rapidly decreasing at infinity together with their derivatives (even Schwartz test functions). The dual space 
${\mathcal S}'_e$ is the space of even Schwartz distributions. Let $H={\mathcal S}'_e$. The partial order in ${\mathcal S}_e$ is given by pointwise comparison of functions, i.e.
$\phi\le \psi$ if $\phi(x)\le \psi(x)$ for every $x\in \R^n.$
\smallbreak
For a function $\phi\in {\mathcal S}_e$ we define the Fourier transform of 
$\phi$ by
$${\mathcal F}\phi(x)=\widehat\phi(x)= \int_{\R^n} \phi(\xi) e^{-i(x,\xi)}\ d\xi,\qquad \forall x\in \R^n,$$
Then ${\mathcal F}: {\mathcal S}_e\to {\mathcal S}_e$ is a bijection.
\smallbreak
We define the {Fourier transform of a distribution} $h\in {{\mathcal S}_e}' $ as a distribution $\widehat{h}$ acting by
$$\langle \widehat{h},\phi \rangle =  \langle h,\widehat\phi \rangle$$
for every $\phi\in {\mathcal S}_e$. Then for every $h\in {\mathcal S}'_e$ and $\phi\in {\mathcal S}_e$
$$
\langle \widehat{h},\widehat\phi\rangle = (2\pi)^n \langle h,\phi \rangle\ {\rm and}\ (\widehat{h})^\wedge= (2\pi)^n h.
$$
An even distribution $h\in {\mathcal S}'_e$ is {\bf positive definite} if its Fourier transform is a positive distribution on even test functions, i.e. 
$\langle {\widehat h}, \phi\rangle \ge 0$ for every non-negative test function $\phi\in {\mathcal S}_e.$  

We say that a distribution is positive on a set $D\subset \R^n$ if it takes positive values on all positive test functions with support in $D.$
We refer the reader to \cite{GS, Ru} for details about Schwartz distributions.

\begin{pr} ${\rm Pos}(\mathcal F)={\rm Pos}({\mathcal{F}}; {\mathcal{S}}_e, {\mathcal{S}}_e, {{\mathcal S}_e}') $  is the set of even positive definite distributions.
\end{pr}
\begin{pf} Suppose that $h\in {\rm Pos}(\mathcal F).$ Let $\phi\in {\mathcal S}_e$ be a non-negative even test function. Since
${\mathcal F}(\widehat{\phi})=(2\pi)^n \phi\ge 0,$ we have $\langle h, \widehat{\phi}\rangle \ge 0.$ This means that $\langle \widehat{h}, \phi\rangle \ge 0.$ Since $\phi$ is an arbitrary non-negative even test function, we see that $\widehat{h}$ is a positive distribution, so $h$ is positive definite.
\smallbreak
On the other hand, suppose that $h$ is an even positive definite distribution. Then for every even test function $\phi$ such that ${\mathcal F}\phi = \widehat{\phi} \ge 0,$ we have
$\langle \widehat{h},\widehat{\phi}\rangle \ge 0,$ which implies $\langle h,\phi \rangle \ge 0.$ This means that $h\in {\rm Pos}(\mathcal F).$
\end{pf}

Positively associated functions have appeared implicitly in convex geometry; see \cite{Ga1, GLW, SW, W1, W2, Z1}. In particular, Goodey, Lutwak and Weil proved in \cite{GLW} that positive functions $h$ on the sphere satisfying $\int_{S^{n-1}} f(x)h(x) dx\ge 0$ for every $f$ with positive spherical Radon transform are radial functions of intersection bodies.
They also proved that a convex body $K$ in $\R^n$ is a projection body if and only if the support function $h_K$ of $K$ satisfies the inequality 
$\int_{S^{n-1}} h_K(x) d\mu(x)\ge 0$ for every even signed measure $\mu$ on the sphere such that the cosine transform $C$ of $\mu$ is positive, i.e.
$$C\mu(x)=\int_{S^{n-1}} |(x,\xi)| d\mu(\xi)\ge 0,\quad \forall x\in S^{n-1}.$$ The authors of \cite{GLW} used these results to characterize the classes of intersection and projection bodies in terms of dual mixed volumes and mixed volumes.
\smallbreak
We work in the setting where $X=Y=H=C^{\infty}_e(S^{n-1})$. The partial order in $C_e^\infty(S^{n-1})$ is given by $\phi\le \psi$ if $\phi(x)\le \psi(x)$ for every $x\in S^{n-1}.$ Applying Definition \ref{positive} to this case, we see that a function $h\in C^{\infty}_{e}(G)$ is positively associated with an operator $T: C^{\infty}_e(S^{n-1})\to C^{\infty}_e(S^{n-1})$ if for any $f\in C^{\infty}_{e}(G)$ satisfying $Tf(x)\ge 0$ for every $x\in S^{n-1}$ we have 
$$\int_{S^{n-1}} f(x) h(x) dx \ge 0.$$
Throughout the paper, we say that a function $f$ is positive on a set $D$ if $f(x)\ge 0$ for every $x\in D.$ If $f(x)>0$ for every $x\in D,$ we say that $f$ is strictly positive on $D.$
\smallbreak
We consider the following class of operators. Let $f\in C_e^{\infty}(S^{n-1}).$ Suppose that $q>-n,\ q\neq 0,2,4,...$ and extend $f$ to an even homogeneous of degree $-n-q$ function on $\R^n\setminus \{0\}:$
$$(f\cdot r^{-n-q})(x)= |x|_2^{-n-q}f(x/|x|_2),$$
where $|x|_2$ is the Euclidean norm on $\R^n.$
Then the Fourier transform of the distribution $f\cdot r^{-n-q}$ is equal to $g\cdot r^q$ for some $g\in C_e^{\infty}(S^{n-1});$ see for example \cite[Lemma 3.16]{K1}.
\smallbreak
We define an operator $T_q: C_e^{\infty}(S^{n-1})\to C_e^{\infty}(S^{n-1})$ by 
$$T_qf(\xi)=g(\xi)= (f\cdot r^{-n-q})^\wedge(\xi),\quad \forall f\in C_e^\infty(S^{n-1}),\quad \forall \xi\in S^{n-1}.$$  
This operator is a bijection of $C_e^{\infty}(S^{n-1})$. 

The operators $T_q$ are related to integral transforms that are well-known in convex geometry. 
A simple argument shows that for $q>0,$ $q\neq 2,4,6...$, the operator $T_q$ is the $q$-cosine transform, up to a constant. This argument can be found in \cite{Se,K-1992} and in \cite[Corollary 3.15]{K1}.

\begin{pr} For any $q>0,\ q\neq 2k,\ k\in \N,$ and any $f\in C_e^\infty(S^{n-1}),$
\begin{equation}\label{q-cosine-fourier}
T_q f(\xi) = (f\cdot r^{-n-q})^\wedge (\xi) = \frac {\sqrt{\pi}\Gamma(-q/2)}{2^{q+1}\Gamma((q+1)/2)} C_qf(\xi).
\end{equation}
\end{pr}
Note that the proof in \cite[Corollary 3.15]{K1} establishes that the Fourier transform of $C_qf\cdot r^q$ coincides with $cf\cdot r^{-n-q}$ outside of the origin in $\R^n,$ where $c$ is a constant. This implies that the functions $C_qf\cdot r^q$ and $c(f\cdot r^{-n-q})^\wedge$ coincide up to a polynomial. Since the functions are even, this polynomial is homogeneous of even degree. But the difference $C_qf\cdot r^q - c(f\cdot r^{-n-q})^\wedge$ is homogeneous of degree $q$ where $q$ is not an even integer. This means that the polynomial is equal to zero and the functions $C_qf\cdot r^q$ and $c(f\cdot r^{-n-q})^\wedge$ are equal. The coefficient in (\ref{q-cosine-fourier}) is written in the form of \cite[Lemma 2.18]{K1}.
\smallbreak
In the case where $q=-1,$ we get the spherical Radon transform.
This fact was proved in \cite{K2} (see also \cite[Lemma 3.7]{K1}).

\begin{pr} For any $f\in C(S^{n-1}),$
 $$T_{-1}f= (f\cdot r^{-n+1})^\wedge =\pi Rf\cdot r^{-1}.$$
 \end{pr}
 
\smallbreak
We use the spherical Parseval's formula: for $-n<q,\ q\neq 2,4,6,...$ and any $f,h\in C_e^{\infty}(S^{n-1}),$
$$\int_{S^{n-1}} (f\cdot r^{-n-q})^\wedge(\xi) (h\cdot r^q)^\wedge(\xi) d\xi =
(2\pi)^n \int_{S^{n-1}} f(x)h(x) dx.$$
This formula was introduced in \cite{K8}; see also \cite[Lemma 3.22]{K1}, \cite[Lemma 5.10]{K1}, \cite{GYY1, GYY2}. A shorter proof can be found in \cite{Mi}. This proof is presented in \cite{Mi} for $q\in (-n,0)$, but it remains the same for $q>0,\ q\neq 2,4,6,...$
\smallbreak
We are ready to characterize the sets ${\rm Pos}(T_q).$ We start with the case $q>0,\ q\neq 2,4,6,...$, and, by (\ref{q-cosine-fourier}), we can formulate the theorem for $C_q.$
\begin{theorem} \label{q>0} For $q>0,\ q\neq 2k,\ k\in \N,$ the class ${\rm Pos}(C_q)$ consists of functions $h\in C^{\infty}_e(S^{n-1})$ such that the Fourier transform of the distribution $\frac 1{\Gamma(-q/2)} h\cdot r^q$ is positive on $\R^n\setminus \{0\}.$
\end{theorem}

\begin{pf} Let $h\in C^\infty_e(S^{n-1})$ be such that the Fourier transform of the distribution $\frac 1{\Gamma(-q/2)} h\cdot r^q$ is positive on $\R^n\setminus \{0\}.$ By \cite[Lemma 3.16]{K1}, $(h\cdot r^{q})^\wedge$ is an infinitely smooth positive function on the sphere extended to an even homogeneous function of degree $-n-q$ on $\R^n\setminus \{0\}.$

Let $f\in C_e^{\infty}(S^{n-1})$ be such that $C_qf\ge 0$ everyhwere on $S^{n-1}.$ By (\ref{q-cosine-fourier}),
$$0\le \frac 1{\Gamma(-q/2)} \int_{S^{n-1}} C_q f(\xi) (h\cdot r^q)^\wedge (\xi) d\xi$$
$$= \frac {2^{q+1}\Gamma((q+1)/2)}{\sqrt{\pi}\Gamma(-q/2)}\frac 1{\Gamma(-q/2)}\int_{S^{n-1}} (f\cdot r^{-n-q})^\wedge (\xi)  (h\cdot r^q)^\wedge (\xi) d\xi,$$
and by Parseval's formula
$\int_{S^{n-1}} f(x) h(x) dx \ge 0.$ So, $h\in {\rm Pos}(C_q).$
\smallbreak
On the other hand, suppose that the function $\frac 1{\Gamma(-q/2)} (h\cdot r^q)^\wedge$ is not positive in $\R^n\setminus \{0\}.$  Then there exists an open symmetric set $\Omega\subset S^{n-1}$ on which the function $\frac 1{\Gamma(-q/2)} (h\cdot r^q)^\wedge$ is strictly negative. Consider a function $\phi\in C_{e}^\infty(S^{n-1})$ such that $\phi\ge 0$ on the whole sphere, $\phi$ is not identically zero, and $\phi$ is supported in $\Omega.$ Define the function $\psi$ by
$\frac 1{\Gamma(-q/2)}(\phi\cdot r^q)^\wedge= \psi\cdot r^{-n-q},$  where $\psi\in C_{e}^\infty(S^{n-1})$ (again by \cite[Lemma 3.16]{K1}). 
For every $\theta\in S^{n-1},$ we have
$$C_q\psi(\theta)=\frac {2^{q+1}\Gamma((q+1)/2)}{\sqrt{\pi}\Gamma(-q/2)} (\psi\cdot r^{-n-q})^\wedge(\theta)$$
$$= (2\pi)^n \frac {2^{q+1}\Gamma((q+1)/2)}{\sqrt{\pi}(\Gamma(-q/2))^2}\phi (\theta)\ge 0.$$
On the other hand, by Parseval's formula,
$$(2\pi)^n \int_{S^{n-1}} \psi(\theta) h(\theta) d\theta = \int_{S^{n-1}} (h\cdot r^q)^\wedge(\xi) (\psi\cdot r^{-n-q})^\wedge(\xi) d\xi$$
$$= (2\pi)^{n}\int_{S^{n-1}} \frac 1{\Gamma(-q/2)} (h\cdot r^q)^\wedge(\xi) \phi(\xi) d\xi < 0,$$
because the positive function $\phi$ can be strictly positive only where $\frac 1{\Gamma(-q/2)}(h\cdot r^q)^\wedge$ is strictly negative.
So, $h\notin {\rm Pos}(C_q),$ which completes the proof.
\end{pf}
Let us pass to the case where $q\in (-n,0).$ In this case $\Gamma(-q/2)$ is always positive, so the result looks simpler.
\begin{theorem} \label{q<0} For $-n<q<0,$ the class ${\rm Pos}(T_q)$ consists of  functions $h\in C^{\infty}_e(S^{n-1})$ such that the distribution $h\cdot r^q$ is positive definite.
\end{theorem}

\begin{pf} Suppose $h\in C^\infty_e(S^{n-1})$ and $h\cdot r^{q}$ is a positive definite distribution. By \cite[Lemma 3.16]{K1}, the Fourier transform $(h\cdot r^q)^\wedge$ is a positive $C^\infty_e$-function on the sphere extended to $\R^n\setminus \{0\}$ as an even homogeneous function of degree $-n-q.$ Let $f\in C^\infty_e(S^{n-1})$ be such that $T_qf$ is a positive function on the sphere, so 
$(f\cdot r^{-n-q})^\wedge\ge 0$ everywhere on $S^{n-1}.$
By Parseval's formula on the sphere, we get
 $$ \int_{S^{n-1}} f(\theta) h(\theta) d\theta = \frac 1{(2\pi)^n} \int_{S^{n-1}} (f\cdot r^{-n-q})^\wedge (\xi) (h\cdot r^q)^\wedge(\xi) d\xi \ge 0,$$
so $h\in {\rm Pos}(T_q).$

Conversely, suppose that $h\in C_e^\infty(S^{n-1})$ and $h\cdot r^q$ is not positive definite. By \cite[Lemma 3.16]{K1}, the Fourier transform of the distribution $h\cdot r^{q}$ is a $C^\infty$-function on the sphere extended
to an even homogeneous of degree $-n-q$ function on $\R^n.$ Since this function cannot be positive, there exists an open 
symmetric set $\Omega\subset S^{n-1}$ on which the function $(h\cdot r^{q})^\wedge$ is strictly negative. Consider a function 
$\phi\in C_{e}^\infty(S^{n-1})$ such that $\phi\ge 0$ on the whole sphere and is not identically zero, and such that $\phi$ is supported in $\Omega.$ Again, by \cite[Lemma 3.16]{K1}, the Fourier transform $(\phi\cdot r^q)^\wedge= \psi\cdot r^{-n-q},$ where $\psi\in C_{e}^\infty(S^{n-1}).$ 
Since $\phi$ is an even function, we have $(\psi\cdot r^{-n-q})^\wedge= (2\pi)^n \phi\cdot r^{q}.$

For every $\theta\in S^{n-1},$
$$T_q\psi(\theta)= (\psi\cdot r^{-n-q})^\wedge(\theta)=(2\pi)^{n} (\phi\cdot r^{q}) (\theta)= (2\pi)^{n}\phi(\theta)\ge 0.$$
On the other hand, by Parseval's formula,
$$(2\pi)^n \int_{S^{n-1}} \psi(\theta) h(\theta) d\theta = \int_{S^{n-1}} (h\cdot r^{q})^\wedge(\xi) (\psi\cdot r^{-n-q})^\wedge(\xi) d\xi$$
$$= (2\pi)^{n}\int_{S^{n-1}} (h\cdot r^{q})^\wedge(\xi) \phi(\xi) d\xi < 0,$$
because the positive function $\phi$ can be strictly positive only in $\Omega$ where $(h\cdot r^{q})^\wedge$ is strictly negative.
So, $h\notin {\rm Pos}(T_q).$
\end{pf}

\section{Comparison problems}
In this section we show how the solution to Problem \ref{general} follows from our results on positively associated functions, as claimed in the introduction. We begin with two lemmas expressing this connection. We use the notation of Problem \ref{general}.
In particular, $T$ is an arbitrary linear operator. The results and proofs of Lemmas \ref{pos-dir} and \ref{neg-dir} hold true if we replace $C_e^\infty(S^{n-1})$ by $C_e(S^{n-1})$, and also if we replace the sphere $S^{n-1}$ by an arbitrary compact set in $\R^n.$

\begin{lemma} \label{pos-dir} Let $p>1.$ Suppose that a positive function $f\in C_e^\infty(S^{n-1})$ has the property that  $f^{p-1}\in {\rm Pos}(T).$  If $g\in C_e^\infty(S^{n-1})$ is any positive function satisfying $Tf(x)\le Tg(x)$ for every $x\in S^{n-1},$ then 
$$\|f\|_{L_p(S^{n-1})} \le \|g\|_{L_p(S^{n-1})}.$$
\end{lemma} 

\proof Since $f^{p-1}\in {\rm Pos}(T)$ and $T(g-f)\ge 0$ pointwise in $S^{n-1}$, we have 
$$\int_{S^{n-1}} (g(x)-f(x))f^{p-1}(x)dx \ge 0.$$
Rearranging the latter inequality and applying H\"older's inequality, we get
$$\|f\|_{L_p(S^{n-1})}^p = \int_{S^{n-1}} f^p(x) dx \le \int_{S^{n-1}} f^{p-1}(x) g(x) dx$$
$$\le \left(\int_{S^{n-1}} f^p(x)dx \right)^{\frac{p-1}p}\left(\int_{S^{n-1}} g^p(x)dx \right)^{\frac 1p} =\|f\|_{L_p(S^{n-1})}^{p-1}
\|g\|_{L_p(S^{n-1})}. \qed$$

\bigbreak
\begin{lemma}\label{neg-dir} Let $p>1,$ and suppose that a function $g\in C_e^\infty(S^{n-1})$ is strictly positive and $g^{p-1}\notin {\rm Pos}(T).$ Then there exists a strictly positive function $f\in C_e^\infty(S^{n-1})$ such that $Tf(x)\le Tg(x)$ for every $x\in S^{n-1},$ but
$$\|f\|_{L_p(S^{n-1})} > \|g\|_{L_p(S^{n-1})}.$$
\end{lemma}

\proof Since $g^{p-1}\notin {\rm Pos}(T),$ there exists a function $h\in C_e^\infty(S^{n-1})$ such that $Th\ge 0$ pointwise in $S^{n-1}$ and
\begin{equation}\label{counter}
\int_{S^{n-1}} g^{p-1}(x) h(x) dx < 0.
\end{equation}
Since $g$ is strictly positive, both $g$ and $h$ are continuous, and $S^{n-1}$ is compact, there exists $\e>0$ so that the function $f=g-\e h$ 
is strictly positive on $S^{n-1}.$ 

Since $Th\ge 0$ pointwise in $S^{n-1}$, we have $Tf=Tg-\e Th\le Tg$ pointwise in $S^{n-1}.$
On the other hand, by the definition of $f$ and by inequality (\ref{counter}), 
$$\int_{S^{n-1}} g^{p-1}(x) f(x) dx$$$$=\int_{S^{n-1}} g^{p}(x) dx -\e \int_{S^{n-1}} g^{p-1}(x) h(x) dx> \|g\|^p_{L_p(S^{n-1})}.$$
By H\"older's inequality,
$$\|g\|_{L_p(S^{n-1}}^{p-1} \|f\|_{L_pS^{n-1})}\ge \int_{S^{n-1}} g^{p-1}(x) f(x) dx > \|g\|^p_{L_p(S^{n-1})}. \qed$$
\smallbreak

Lemmas \ref{pos-dir} and \ref{neg-dir} in conjunction with Theorems  \ref{q>0} and \ref{q<0} immediately imply the solution of Problem \ref{general} for the operators $T_q.$ In the case where $q>0,\ q\neq 2,4,6,...$ we get Theorem \ref{comp-q>0} formulated in the introduction.
\smallbreak
In the case $-n<q<0,$ we get the following result.
\begin{theorem}\label{comp-q<0}  Let $p>1,\ -n<q<0.$ 
\smallbreak
(i) Suppose that $f\in C_e^\infty(S^{n-1})$ is a positive function so that $f^{p-1}\cdot r^{q}$ is a positive definite distribution. If $g\in C_e^\infty(S^{n-1})$ is a positive function such that $T_qf\le T_qg$ pointwise on $S^{n-1}$, then
$$\|f\|_{L_p(S^{n-1})} \le \|g\|_{L_p(S^{n-1})}.$$

(ii) If $g\in C^\infty_e(S^{n-1})$ is a strictly positive function such that $g^{p-1}\cdot r^{q}$ is not a positive definite distribution, then there exists a strictly positive function $f\in C^\infty_e(S^{n-1})$ such that $T_qf\le T_qg$ pointwise on $S^{n-1}$, but
$$\|f\|_{L_p(S^{n-1})} > \|g\|_{L_p(S^{n-1})}.$$
\end{theorem}

In the case $q=-1,$ where $\frac 1{\pi} T_{-1}=R$ is the spherical Radon transform, Theorem \ref{comp-q<0}  was proved in [KRZ]
without using positively associated functions.
\medbreak
\noindent{\bf Remark 1.} Note that when $p=1$ in Theorems \ref{comp-q>0} and \ref{comp-q<0}, we get an affirmative answer to Problem 1. In fact, 
if $q>0,$ integrating the inequality $C_qf(\xi)\le C_qg(\xi)$ by $d\xi$ over the sphere, we get
$\int_{S^{n-1}} f(x) dx \le \int_{S^{n-1}} g(x) dx$. 

In the case where $q\in (-n,0),$ note that the Fourier transform of powers of the Euclidean norm is a constant on the sphere; see \cite{GS}. Namely, for every $\xi\in S^{n-1},$
\begin{equation}\label{eucl}
(|x|_2^q)^\wedge(\xi) = \frac{2^{n+q}\pi^{n/2}\Gamma(\frac {n+q}2)}{\Gamma(-\frac q2)}:= \frac 1{d_{n,q}}.
\end{equation}
Now, applying Parseval's formula we get
$$\int_{S^{n-1}} T_q f(x) dx= d_{n,q} \int_{S^{n-1}} (f\cdot r^{-n-q})^\wedge(\xi) (|x|_2^{q})^\wedge(\xi) d\xi$$$$=
d_{n,q}(2\pi)^n \int_{S^{n-1}} f(x) dx.$$
It is clear now that the inequalities $T_qf(x)\le T_qg(x)$ for all $x\in S^{n-1}$ imply $\int_{S^{n-1}} f(x) dx \le \int_{S^{n-1}} g(x) dx.$
\medbreak

\noindent{\bf Remark 2.} Let ${\bf 1}$ be the function identically equal to 1 on $S^{n-1}.$  Lemma \ref{pos-dir} implies a lower estimate for an arbitrary linear operator $T$ acting from $C_e^\infty(S^{n-1})$ to itself. We call this estimate the slicing inequality for the operator $T.$
 
\begin{co}\label{slicing-T}  Suppose that $T{\bf 1}= c_T{\bf 1},$ where $c_T>0.$ 
Let $p>1,$ and let $f\in C_e^\infty(S^{n-1})$ be a positive function such that $f^{p-1}\in {\rm Pos}(T).$
Then
\begin{equation} \label{lower-T}
\|f\|_{L_p(S^{n-1})}\le \frac{|S^{n-1}|^{\frac 1p}}{c_T} \max_{z\in S^{n-1}} Tf(z).
\end{equation}
\end{co}

\proof First, $\max_{z\in S^{n-1}} Tf(z) > 0$ unless $f\equiv 0.$ In fact, if $Tf(z)\le 0$ for every $z\in S^{n-1},$ then, since $f^{p-1}\in {\rm Pos}(T),$ we have $\int_{S^{n-1}} f^{p}(x) dx\le 0$. But $f$ is a positive continuous function on $S^{n-1}$, so $f\equiv 0.$

Let $g$ be a constant function on $S^{n-1}$ defined by
$$g=\frac 1{c_T} \max_{z\in S^{n-1}} Tf(z){\bf 1}.$$
Then $g$ is positive on $S^{n-1}$ and $Tg(x)= \max_{z\in S^{n-1}} Tf(z) \ge Tf(x)$ for every $x\in S^{n-1}.$ By Lemma \ref{pos-dir},
$$\|f\|_{L_p(S^{n-1})}\le \|g\|_{L_p(S^{n-1})}= \frac{|S^{n-1}|^{\frac 1p}}{c_T} \max_{z\in S^{n-1}} Tf(z). \qed$$

For the operators $T_q,$ we have
$$\frac 1{\Gamma(-q/2)} T_q{\bf 1}(\xi)= \frac 1{\Gamma(-q/2)}({\bf 1}\cdot r^{-n-q})^\wedge(\xi)$$$$= \frac 1{\Gamma(-q/2)}(|x|_2^{-n-q})^\wedge(\xi)= \frac{2^{-q}\pi^{\frac n2}}{\Gamma(\frac{n+q}2)}{\bf 1}:= c_{n,q}{\bf 1}.$$
for every $\xi\in S^{n-1}.$ Here we used formula (\ref{eucl}). The latter equality and Corollary \ref{slicing-T} imply a slicing inequality for the operators $T_q,\ q>-n,\ q\neq 0,2,4,...$
\medbreak

\noindent{\bf Remark 3.} For strictly positive functions $f$ the condition of Theorem \ref{comp-q>0} that the Fourier transform of  
$\frac 1{\Gamma(-q/2)} f^{p-1}\cdot r^q$ is a positive distribution on $\R^n\setminus \{0\}$ is equivalent to the function $f$ being equal to $f=\|\cdot\|^{\frac q{p-1}}, $ where the norm belongs to an $n$-dimensional subspace of $L_q(S^{n-1}),\ q>0,\ q\neq 2,4,6,...$; see \cite{K-1992} or \cite[Theorem 6.10]{K1}. The condition of Theorem \ref{comp-q<0} that $f^{p-1}\cdot r^q$ is a positive definite distribution on $\R^n$ is equivalent to $f=\|\cdot\|^{\frac q{p-1}},$ where the norm belongs to an $n$-dimensional space that embeds in $L_q,\ q<0;$ see \cite[Section 6.3]{K1}. Many examples of functions $f$ satisfying and not satisfying these conditions can be found in \cite[Chapter 6]{K1}. However, this connection seems to deserve more attention and will be considered in a future paper.

\end{document}